\documentclass[1p]{elsarticle}
\usepackage{amssymb}
\usepackage{amsfonts}

\newtheorem{theorem}{Theorem}

\newtheorem{lemma}[theorem]{Lemma}

\newproof{pf}{Proof}

\begin{document}

\title{On triangle-free list assignments}

\author[agh]{Jakub Przyby{\l}o}

\address[agh]{AGH University of Science and Technology, Faculty of Applied Mathematics, al. A. Mickiewicza 30, 30-059 Krakow, Poland}

\begin{abstract}
We show that Bernshteyn's proof of the breakthrough result of Molloy that triangle-free graphs are choosable from lists of size $(1+o(1))\Delta/\log\Delta$ can be adapted to yield a stronger result. In particular one may prove that such list sizes are sufficient to colour any graph of maximum degree $\Delta$ provided that  vertices sharing a common colour in their lists do not induce a triangle in $G$, which encompasses all cases covered by Molloy's theorem. This was thus far known to be true for lists of size $(1000+o(1))\Delta/\log\Delta$, as implies a more general result due to Amini and Reed. We also prove that lists of length $2(r-2)\Delta \log_2\log_2\Delta/\log_2\Delta$ are sufficient if one replaces the triangle by any $K_r$ with $r\geq 4$, pushing also slightly the multiplicative factor of $200r$ from Bernshteyn's result down to $2(r-2)$. All bounds presented are also valid within the more general setting of correspondence colourings.
\end{abstract}

\begin{keyword}
triangle-free graphs \sep list colouring \sep DP-colouring \sep triangle-free list assignment \sep Johansson's Theorem
\end{keyword}

\maketitle

\section{Introduction}

The chromatic number of any graph $G$ is trivially bounded above by $\Delta+1$, where $\Delta$ is the maximum degree of $G$. This bound is tight in general. It was presumably Vizing~\cite{Vizing68} who first asked if it can be improved in case of triangle-free graphs.
The first non-trivial results confirming this supposition were due to Borodin and Kostochka~\cite{BorodinKostochka-triangle-free}, Catlin~\cite{Catlin-triangle-free}, and Lawrence~\cite{Lawrence-triangle-free}, who independenty proved that $\chi(G)\leq \frac{3}{4}(\Delta+2)$ for such graphs, which was later improved to $\chi(G)\leq \frac{2}{3}\Delta+2$ by Kostochka~\cite{Jensen-Toft-book}.

The quest for a proper colouring of a triangle-free graph relates intuitively
with the well-known Coupon Collector's Problem~\cite{SedgewickFlajolet}, which
asks how many independent draws from a fixed set of coupons, say $[n]:=\{1,2,\ldots,n\}$ assure 
collecting the whole deck $[n]$ with high probability.  It is known that the threshold in this problem is roughly $m=n\log n$ draws. In particular, if $m=(1-\varepsilon)n\log n$, where `$\log$' denotes the natural logarithm, then w.h.p. one shall be relatively far from collecting all coupons.
Knowing this imagine we are randomly choosing one of $n$ colours for each of the neighbours of a given vertex $v$ of degree $m=(1-\varepsilon)n\log n$, disregarding for now the rest of our triangle-free graph. No colour conflict shall arise this way, as the neighbourhood induces an edge-less graph. Moreover w.h.p. there shall still be an available choice of colour, in fact many of these, for $v$. 
This observation suggested a possibly achievable value of an upper bound for the chromatic number of a triangle-free graph with maximum degree $\Delta=m$, that is roughly $\Delta/\log \Delta$ (as $n\log n \approx\Delta$ for $n\approx \Delta/\log \Delta$). Many efforts were devoted to finally turn this heuristic intuition into a proof of an $O(\Delta/\log \Delta)$ upper bound by Johansson~\cite{Johansson-triangle-free}, who improved and extended an earlier result of Kim~\cite{Kim} encompassing exclusively graphs of girth at least $5$.
Johansson provided a technically elaborate probabilistic argument, based on iterative application of the Lov\'asz Local Lemma, which in fact implied that also in a more general list setting: $\chi_l(G)\leq 9\Delta/\log \Delta$ for every triangle-free graph $G$, which was later refined to $\chi_l(G)\leq 4\Delta/\log \Delta$ in~\cite{PettieSu}.

One of the key issues which had to be resolved by Johansson  was the fact that lists analyzed within his multistage random process could differ to a large extent, and more importantly their evolution was not independent -- in particular dependency between lists of colours available for any two naighbours of a given vertex  could potentially be much stronger than in the case of graphs of girth at least $5$  (though intuitively a large common neighbourhood of such two neigbours of a given vertex $v$ should increase the chance for assigning the same colour to these two neighbours, and consequently should result in less necessary colours in the list of $v$, which nevertheless provided more technical problems than could have been until recently beneficially utilized). 
Many of these issues were lately very efficiently overcome by Molloy~\cite{Molloy-small_clique_number} within his remarkable and concise proof, based on the \emph{entropy compression} method (name due to Tao~\cite{Tao}), developed in the wake of Moser and Tardos study over algorithmic version of the Lov\'asz Local Lemma~\cite{Moser,MoserTardos} and later variations and modifications, see e.g.~\cite{GrytczukKozikMicek,AchlioptasIliopoulos,BosekCzerwinskiGrytczukRzazewski,DujmoviJoretKozikWood,EsperetParreau}. As a result Molloy proved that $\chi_l(G)\leq (1+o(1))\Delta/\log \Delta$, drawing very near to the lower bound $(1/2-o(1))\Delta/\log \Delta$ which derives from investigating random $\Delta$-regular graphs (applicable not only to triangle-free graphs, but also graphs with girth larger than any fixed constant $k\geq 3$), cf.~\cite{Bollobas-chromatic,Kim,Molloy-small_clique_number}.
Molloy's approach was also very recently turned into even simpler and shorter argument  by Bernshteyn~\cite{Bernshteyn}, who used a very inventive and innovative application of the Lov\'asz Local Lemma to provide a one-round (i.e without iterations) proof of the results from~\cite{Molloy-small_clique_number}, extending these at the same time towards a more general setting of DP-colouring (known also as \emph{correspondence coloring}), introduced quite recently by Dvo\v{r}\'ak and Postle~\cite{DvorakPostle} in order to settle a long-standing open question
of Borodin regarding list colouring of planar graphs with no cycles of certain lengths~\cite{Borodin-plane}, see also~\cite{MolloyPostle2021} for most recent results and up-to-date exposition of this subject.

One of the key ideas facilitating Molloy's (and Bernshteyn's) approach is
a new result on a variation of the Coupon Collector's Problem, which we formalize below,
implying in particular that a similar threshold (under which one is very likely not to collect a full deck) as discussed above is also valid in case of drawing coupons from subdecks of possibly highly diversified sizes as long as each of these has at least two elements. 
This last requirement fits in however perfectly in Molloy's idea consisting in introducing 
blank lottery draws (which correspond to choosing a partial colouring, where vertices associated with blanks have no coulour assigned). 

Both, Molloy~\cite{Molloy-small_clique_number}  and Bernshteyn~\cite{Bernshteyn} investigated
$K_r$-free graphs with $r$ bigger than $3$ as well. They obtained upper bounds of order larger by a factor $\log\log\Delta$ than the suspected by Alon, Krivelevich and Sudakov~\cite{AlonKrivelevichSudakov} optimal one. In particular Bernshteyn proved that for any fixed $r\geq 4$, $\chi_{DP}(G)\leq 200r\Delta \log_2\log_2\Delta/\log_2\Delta$ for each $K_r$-free graph $G$ with maximum degree large enough.

In the present paper we show that almost exactly the same technique as applied by Bernshteyn~\cite{Bernshteyn} yields a stronger result.
Namely we observe that lists of size $(1+o(1))\Delta/\log\Delta$ are sufficient not only for triangle-free graphs, but also for any other graph $G$ as long as vertices sharing a common colour
in their lists do not induce a triangle in $G$ (this obviously encompasses all cases covered by Molloy's result~\cite{Molloy-small_clique_number}), cf. Theorem~\ref{Triangl-Free-Epsilon-Theorem} below. We moreover show that the same holds in case of more general DP-colourings, cf. Theorem~\ref{Triangl-Free-DP-Theorem}, which includes all cases from Bernshteyn's result~\cite{Bernshteyn} as well. Finally we also similarly extend Molloy's and Bernshteyn's results for $K_r$-free graphs, pushing at the same the factor $200r$ down to $2(r-2)$, cf.~Theorem~\ref{Kr-Free-DP-Theorem}.

Somewhat related direction of research can be found in~\cite{Reed-list_const}, where Reed
turned towards transferring certain requirements guaranteeing existence of a proper vertex colouring from graphs to list arrangements. Namely, he conjectured that there
exists a proper colouring of a graph $G$ from any assignment of lists of size larger than $d$, provided that for each vertex $v$ and any colour $c$ in tits list, the number of neighbours of $v$ whose lists also contain $c$ is at most $d$, regardless of the maximum degree of $G$ itself,
and showed it is true for lists of size at least $Cd$, where $C=2e\approx 5.44$.
This conjecture in its verbatim meaning was disproved by Bohman and Holzman~\cite{BohmanHolzman}, and at the same time confirmed in the asymptotic sense
by Reed and Sudakov~\cite{ReedSudakov} (via improving $C$ to $1+o(1)$).
 Alon and Assadi~\cite{AlonAssadi} in turn proved that lists of size $(C'+o(1))d/\log d$ with $C'=8$ are sufficient for triangle-free graphs. 
 The constant $C'$ was further pushed down to $1$ by Cambie and Kang, but only for the restricted bipartite case, cf.~\cite{CambieKang}.
Under additional stronger assumption, the same as investigated in the present paper, i.e. that vertices sharing a common colour in their lists do not induce a triangle in $G$, Amini and Reed~\cite{AminiReed} proved 
that lists of size $(1000+o(1))d/\log d$ are sufficient; 
note this implies Theorem~\ref{Triangl-Free-Epsilon-Theorem}, but with a worst constant, i.e. $1000$ instead of $1$.

Yet another interesting result related with lists' structure is an extension of Galvin's theorem confirming the famous List Colouring Conjecture for bipartite graphs due to Fleiner~\cite{Fleiner}, who proves that the conjecture holds in case of any graph and each edge list assignment with no odd cycle having a common colour in the lists of its edges.

In the following section we briefly discuss rather standard tools of the probabilistic method.
Next we formalize a variant of Molloy's generalization of the Coupon Collector's Problem, dissecting it from the main argument as an entity interesting by its own. We believe this brings out its role in the whole argument, makes the reasoning more clear, and emphasises simplicity of Bernshteyn's and Molloy's approach, compared to original Johanson's probabilistic proof. We present all details of this relatively plain argument in the subsequent section in the more general setting of lists (rather than just graphs) inducing no triangles.  
We next discuss how it extends to DP-colourings, and in the final section we present the result for $K_r$-free list assignments directly for DP-colourings.

\section{Preliminaries}

Given a graph $G=(V,E)$, $v\in V$ and $A,B\subseteq V$, we denote by $N_G[v]$ the \emph{closed neigbhbourhood} of $v$, i.e. $N_G(v)\cup\{v\}$, writing $N[v]$ instead if this causes no ambiguities, and we set $N_G(A):=\bigcup_{u\in A}N_G(u)$, $N_G[A]:=\bigcup_{u\in A}N_G[u]$. Moreover, we denote by $E_G[A,B]$ the set of edges joining $A$ with $B$ in $G$, and by $G-A$ the subgraph $G[V\smallsetminus A]$ induced by $V\smallsetminus A$ in $G$.

We say that $\{0,1\}$-valued random variables $Y_1,Y_2,\ldots,Y_n$ are \emph{negatively correlated} if for every $S\subseteq [n]$:
\begin{equation}\label{NegativeCorrelationCondition1}
\mathbf{Pr}\left(\bigwedge_{i\in S} \{Y_i=1\}\right)\leq \prod_{i\in S} \mathbf{Pr}\left(Y_i=1\right).
\end{equation}
Note~(\ref{NegativeCorrelationCondition1})  follows in particular by the condition: for every $j\in [n]$ and each $S\subseteq [n]\smallsetminus\{j\}$:
$$\mathbf{Pr}\left(Y_j=1~|~\bigwedge_{i\in S} \{Y_i=1\}\right)\leq  \mathbf{Pr}\left(Y_j=1\right).$$
As in turn $\mathbf{Pr}(A|B)\leq \mathbf{Pr}(A)$ is equivalent to  
$\mathbf{Pr}(B|\bar{A})\geq \mathbf{Pr}(B)$ for any events $A,B$, negative correlation of $Y_1,Y_2,\ldots,Y_n$ follows if  for every $j\in [n]$ and each $S\subseteq [n]\smallsetminus\{j\}$:
\begin{equation}\label{NegativeCorrelationCondition3}
\mathbf{Pr}\left(\bigwedge_{i\in S} \{Y_i=1\}~|~Y_j=0\right)\geq  \mathbf{Pr}\left(\bigwedge_{i\in S} \{Y_i=1\}\right).
\end{equation}

\begin{lemma}[Chernoff Bounds, \cite{PanconesiSrinivasan}]\label{ChernoffBoundNegativeCorrelation}
Let  $X_1,X_2,\ldots,X_n$ be $\{0,1\}$-valued random variables. Set $X:=\sum_{i=1}^nX_i$ and $Y_i:=1-X_i$ for $i\in [n]$. Then for every $0<a\leq 1$:
\begin{enumerate}
\item[(i)] if $X_1,X_2,\ldots,X_n$  are negatively correlated, then 
$$\mathbf{Pr}\left(X > (1+a)\mathbf{E}(X)\right) < e^{-\frac{a^2\mathbf{E}(X)}{3}};$$
\item[(ii)] if $Y_1,Y_2,\ldots,Y_n$  are negatively correlated, then 
$$\mathbf{Pr}\left(X < (1-a)\mathbf{E}(X)\right) < e^{-\frac{a^2\mathbf{E}(X)}{2}}.$$
\end{enumerate}
\end{lemma}

We shall use the symmetric variant of the Lopsided Lov\'asz Local Lemma, see e.g.~\cite{AlonSpencer} 
(in particular Corollary 5.1.2 and the comments below). 
\begin{theorem}[\textbf{Symmetric Lopsided Local Lemma}]
\label{LLL-symmetric-lopsided}
Let $\mathcal{A}=\{A_1,A_2,\ldots,A_n\}$ be a family of 
events in an arbitrary pro\-ba\-bi\-li\-ty space.
Suppose that for every $i\in [n]$ there is a set $\Gamma(A_i)\subset \mathcal{A}$ of size at most $D$ such that for each $\mathcal{B}\subseteq \mathcal{A}\smallsetminus (\Gamma(A_i)\cup \{A_i\})$: 
$$ \mathbf{Pr}(A_i|\bigcap_{B\in\mathcal{B}}\overline{B})  \leq \frac{1}{e(D+1)}.$$
Then $ \mathbf{Pr}\left(\bigcap_{i=1}^n\overline{A_i}\right)>0$.
\end{theorem}

\section{Molloy's Coupon Collector's Problem}
\label{MolloysCouponCollectorsProblem}

We use integer $0$ to represent the \emph{blank coupon}.
A set $L\subseteq [n]$, $L\neq \emptyset$ shall be called an \emph{$n$-deck}, 
while ${^0L}:=\{0\}\cup L$ -- an \emph{$n$-deck with a blank}.
Suppose that given $n$-decks $L_1,L_2,\ldots,L_m$, we are independently choosing a random element from every 
${^0L_i}$ -- each with equal probability (i.e. $1/|{^0L_i}|=1/(|L_i|+1)$). We call such a process 
an \emph{$(L_1,L_2,\ldots,L_m)$-lottery with blanks}. Let  $c_1,c_2,\ldots,c_m$ be the elements chosen from the corresponding sets within such lottery. 
We say that \emph{$L_i$ is missed} if $c_i=0$. 
(Such notation relates also with a buildup towards a possible second stage of the lottery, during which one reuses the decks with a blank chosen in the first phase, what seems potentially applicable 
within some further exploitations of this process). 
We call $U:=[n]\smallsetminus\{c_1,c_2,\ldots,c_m\}$ the set of \emph{uncollected coupons}.
The proof of the following lemma mimics the arguments in~\cite{Molloy-small_clique_number}  and~\cite{Bernshteyn}.
It states that if $m$ is appropriately smaller than $n\log n$, then w.h.p. many elements of the full deck $[n]$ shall remain uncollected, while neither of these uncollected coupons shall be available in many lists which were missed. 

\begin{lemma}\label{CouponCollectorsWithBlanksLemma}
For every $\varepsilon\in(0,1)$ there exists $n_\varepsilon$ such that for any $n$-decks $L_1,L_2,\ldots,L_m$ with $m\leq (1-\varepsilon)n\log n$ and $n\geq n_\varepsilon$, if $U$ is the set of uncollected coupons within the $(L_1,L_2,\ldots,L_m)$-lottery with blanks, then
\begin{itemize}
\item[(i)] $\mathbf{Pr}\left(|U|<(1-\varepsilon)n^\varepsilon\right)<e^{-\frac{\varepsilon^2n^\varepsilon}{2}}$; 
\item[(ii)] for every $c\in [n]$: $\mathbf{Pr}\left(c\in U ~{\rm and}~ c~{\rm belongs ~to ~more ~than~}\varepsilon^2n^\varepsilon {\rm~missed~lists}\right)<e^{-\frac{\varepsilon^2n^\varepsilon}{6}}$. 
\end{itemize}
\end{lemma}

\begin{pf} Let $\mathcal{L}=\{L_1,L_2,\ldots,L_m\}$.
Consider any $c\in [n]$ and let $\mathcal{L}_c:=\{L\in \mathcal{L}:c\in L\}$ be the set of these given $n$-decks which contain the coupon $c$.
Since the choices within our lottery are independent, then
$$\mathbf{Pr}\left(c\in U\right) = \prod_{L\in\mathcal{L}_c}\frac{|L|}{|L|+1}
= \prod_{L\in\mathcal{L}_c}\left(1 - \frac{1}{|L|+1}\right),$$ 
hence
\begin{equation}\label{cinUestimation}
e^{-\sum_{L\in\mathcal{L}_c}\frac{1}{|L|}}\leq \prod_{L\in\mathcal{L}_c}e^{-\frac{1}{|L|}} \leq\mathbf{Pr}\left(c\in U\right)\leq \prod_{L\in\mathcal{L}_c}e^{-\frac{1}{|L|+1}}\leq e^{-\sum_{L\in\mathcal{L}_c}\frac{1}{|L|+1}}
\end{equation}
(where a product and a sum over the empty set are understood as $1$ and $0$, respectively). Thus,
by convexity of the function $e^{-x}$ for $x>0$ and the fact that $L\subseteq [n]$ for $L\in\mathcal{L}$:
\begin{eqnarray}
\mathbf{E}(|U|)&\geq& \sum_{c\in[n]}e^{-\sum_{L\in\mathcal{L}_c}\frac{1}{|L|}} 
= n \left(\frac1n \sum_{c\in[n]}e^{-\sum_{L\in\mathcal{L}_c}\frac{1}{|L|}} \right)
\geq n e^{-\frac1n \sum_{c\in[n]}\sum_{L\in\mathcal{L}_c}\frac{1}{|L|}} \nonumber\\
&=& n e^{-\frac1n \sum_{L\in\mathcal{L}}\sum_{c\in L}\frac{1}{|L|}} 
= ne^{-\frac{m}{n}}
\geq ne^{-\frac{(1-\varepsilon)n\log n}{n}}
= n^\varepsilon. \label{ExpectationOfU-estimation}
\end{eqnarray}
Let $X_c$ denote the random variable valued $1$ when $c\in U$ and $0$ otherwise. Thus $|U|=\sum_{i=1}^nX_i$. 
We shall show that the random variables $Y_1,Y_2,\ldots,Y_n$ where $Y_i:=1-X_i$ for $i\in[n]$ are negatively correlated.
It is sufficient to observe that for every $j\in [n]$ and each $S\subseteq [n]\smallsetminus\{j\}$ condition (\ref{NegativeCorrelationCondition3}) holds, i.e. that:
$$
\mathbf{Pr}\left(S\cap U=\emptyset~|~j\in U\right)\geq  \mathbf{Pr}\left(S\cap U=\emptyset\right), \nonumber
$$
which is equivalent to the following one, where $c_i$ denotes the element drawn from $L_i$ for every $i\in[m]$:
\begin{equation}\label{NegativeCorrelationCondition3-application_2}
\mathbf{Pr}\left(\forall k\in S ~\exists i\in [m]: ~c_i=k~|~\forall i\in[m]:~c_i\neq j\right)\geq  \mathbf{Pr}\left(\forall k\in S ~\exists i\in [m]: ~c_i=k\right).
\end{equation}
This is rather straightforward, as the fact that $j$ cannot be chosen from any given $n$-deck only increases the probability that other types of coupons are drawn (we in a sense commit choices from shorter lists).
More formally, one may mimic the random process corresponding to the probability on the left-hand side of~(\ref{NegativeCorrelationCondition3-application_2}) 
using the original unconditioned one, where one by one we make independent choices from consecutive complete decks but whenever one chooses $j$ from a given deck $L_i$, then they discard such a coupon and redraw a new element, this time from $L_i\smallsetminus\{j\}$ with a uniform probability distribution. 
Since any result $(c_1,\ldots,c_m)$ of the original process satisfying: $\{\forall k\in S ~\exists i\in [m]: ~c_i=k\}$ must also satisfy it within our modified procedure (as $j\notin S$),  
the probability on the left-hand side of~(\ref{NegativeCorrelationCondition3-application_2}) must be at least as large as the one on the right. 

By Lemma~\ref{ChernoffBoundNegativeCorrelation}, due to the negative correlation of $Y_1,Y_2,\ldots,Y_n$ and~(\ref{ExpectationOfU-estimation}) we thus have: 
$$\mathbf{Pr}\left(|U|<(1-\varepsilon)n^\varepsilon\right)<e^{-\frac{\varepsilon^2n^\varepsilon}{2}},$$
as desired.

In order to show that (ii) holds for any given $c\in[n]$, let us denote by $\mathcal{L'}_c=\{L_i:c\in L_i, c_i=0\}$
the set of missed $n$-decks including $c$. 
We thus wish to estimate the probability of the event 
$$A_c: ~c\in U~\wedge~\left|\mathcal{L'}_c\right|>\varepsilon^2n^\varepsilon.$$
We shall actually show that either it was highly unlikely already from the start that $c\in U$, i.e. that $c$ was not collected, or otherwise there must have been relatively few $n$-decks including $c$ (or possibly more of these but of larger sizes) out of which only a small enough fraction could have been missed w.h.p. 
We may in particular assume that 
$$\sum_{L\in\mathcal{L}_c}\frac{1}{|L|+1}\leq \frac{1}{2}\varepsilon^2n^\varepsilon,$$
as otherwise, by~(\ref{cinUestimation}),
$$\mathbf{Pr}(A_c)\leq \mathbf{Pr}(c\in U) \leq e^{-\sum_{L\in\mathcal{L}_c}\frac{1}{|L|+1}} < e^{-\frac{\varepsilon^2n^\varepsilon}{2}} < e^{-\frac{\varepsilon^2n^\varepsilon}{6}},$$
as desired. We thus have that
$$\mathbf{E}\left(|\mathcal{L'}_c|\right)=\sum_{L_i\in \mathcal{L}_c}\mathbf{Pr}(c_i=0) 
= \sum_{L\in \mathcal{L}_c}\frac{1}{|L|+1}
\leq \frac{1}{2}\varepsilon^2n^\varepsilon.$$
Therefore, by the Chernoff Bound,
$$\mathbf{Pr}\left(A_c\right)\leq \mathbf{Pr}\left(\left|\mathcal{L'}_c\right|>\varepsilon^2n^\varepsilon\right) 
< e^{-\frac{\varepsilon^2n^\varepsilon}{6}}.$$
\qed
\end{pf}

\section{Triangle-free list assignments}

Let $G=(V,E)$ be a graph and consider a \emph{list assignment} $L:V\to 2^{\mathbb{N}_+}$,
associating to every vertex a possibly empty set of positive integers; we call it a \emph{list $k$-assignment} if $|L(v)|=k$ for $v\in V$. 
Set ${^0L}(v):=\{0\}\cup L(v)$ for every $v\in V$.
A \emph{partial $L$-colouring} of $G$ is an assignment $\omega:V\to\mathbb{N}$ such that 
$\omega(v)\in {^0L}(v)$  for every $v\in V$ and $\omega(u)\neq \omega(v)$ for each $uv\in E$ with $\omega(u),\omega(v)\neq 0$.
The vertices assigned $0$ are interpreted as uncoloured, we denote the set of all such vertices by $O_\omega$ and set $O_\omega(u):=O_\omega\cap N(u)$ for $u\in V$. 
If there is such $\omega$ with $\omega(v)>0$ for every $v\in V$, we call it an \emph{$L$-colouring} of $G$ while $G$ is called \emph{$L$-colourable}.

Let $L$ be a list assignment of $G$. 
We say that a colour $c$ is \emph{available} for $v$ if $c\in L(v)$.  
Given a graph $F$, the list assignment $L$ is said to be \emph{$F$-free}
if for every colour $c>0$, $F$ is not a subgraph of 
the graph induced in $G$ by the vertices with colour $c$ available 
(or equivalently, if the vertices of no subgraph of $G$ isomorphic to $F$ have a common available colour).
For any $u\in V$ and $c\in L(u)$, we shall denote by $$N_L^c(u):=\{v\in N(u):c\in L(v)\}$$ the set of neighbours of $u$ with $c$ available for them, while $|N_L^c(u)|$ shall be addressed as \emph{availability of $c$ around $u$}.
Starting from the initial list assignment with colour availabilities potentially much larger than list sizes,  we shall find a partial colouring such that in what remains, availabilities shall be significantly smaller than the sizes of list leftovers (for uncoloured vertices).
This shall occur sufficient e.g. due to the mentioned result of Reed~\cite{Reed-list_const}.
\begin{theorem}[Reed~\cite{Reed-list_const}]\label{Lists8TimesAvailiabilityObs}
A graph $G=(V,E)$ is $L$-colourable for any list assignment $L$ of $G$ such that for some fixed constant $l$, $|L(u)|\geq l$ and $|N_L^c(u)|\leq l/(2e)$ for each $u\in V$ and $c\in L(u)$.
\end{theorem}
A colouring $\omega$ consistent with Theorem~\ref{Lists8TimesAvailiabilityObs} can be found randomly. Once one learns from~\cite{MolloyReed} that it is sufficient to define bad events as 
$A_{uv,c}:~\omega(u)=c=\omega(v)$ for every $uv\in E$ and $c\in L(u)\cap L(v)$, the rest is a simple exercise based on application of the Local Lemma, see e.g.~\cite{MolloyReed,Reed-list_const} for details (where it is convenient to assume a uniform list size). We might have used a few other results mentioned in the introduction, with better upper bounds, but Theorem~\ref{Lists8TimesAvailiabilityObs} is completely sufficient for our purposes.

The following theorem generalizes Molloy's result~\cite{Molloy-small_clique_number} towards triangle-free list assignments with list sizes $(1+o(1))\Delta/\log\Delta$. Its proof is an extension of an approach proposed by Bernshteyn~\cite{Bernshteyn}.   

\begin{theorem}\label{Triangl-Free-Epsilon-Theorem}
For any $0<\varepsilon<1/3$ there exists $\Delta_\varepsilon$ such that every graph $G$ 
with maximum degree $\Delta\geq \Delta_\varepsilon$ is $L$-colourable for any triangle-free list assignment $L$ of $G$ as long as each of the lists has size at least 
$\left\lfloor(1+4\varepsilon)\Delta/\log\Delta\right\rfloor$.
\end{theorem}

\begin{pf}
Fix any $\varepsilon\in (0,1/3)$ and $G=(V,E)$ consistent with the assumption of the theorem; we do not specify $\Delta_\varepsilon$, we only assume that $\Delta$ is large enough so that all explicit inequalities below hold.
For any partial colouring $\omega:V\to\mathbb{N}$ and any list assignment $L$ of $G$ 
we define $L^\omega$ to be the list assignment (of $G$) of still available colours (for vertices uncoloured within $\omega$), i.e.
$$L^\omega(u)=\left\{\begin{array}{lcr}
L(u)\smallsetminus\{\omega(v):v\in N(u)\} & if & \omega(u) = 0\\
\emptyset &if & \omega(u) > 0
\end{array}
\right..$$
Set $n:=\left\lfloor(1+4\varepsilon)\Delta/\log\Delta\right\rfloor$, $l:=(1-\varepsilon)n^\varepsilon$ and assume from now on that $L$ is a list $n$-assignment of $G$, i.e. $|L(v)|=n$ for every $v$.  
Suppose $\omega$ is a partial $L$-colouring of $G$.
Note that if 
\begin{equation}\label{MainGoal}
|L^\omega(u)|\geq l~~ {\rm and}~~ |N_{L^\omega}^c(u)|\leq l/(2e)
\end{equation}
for every $c\in L^\omega(u)$ and each $u\in O_\omega$, then the assertion follows from Theorem~\ref{Lists8TimesAvailiabilityObs}, for it is sufficient to find a proper colouring of the vertices in $O_\omega$ from their lists $L^\omega(u)$, $u\in O_\omega$ 
in order to complete the construction of an $L$-colouring of $G$ 
in such a case.

We may interpret $L$ as a binary relation and define any $L'\subseteq L$ to be an 
\emph{$L$-subassignment}, i.e. a list assignment $L'$ of $G$ such that $L'(v)\subseteq L(v)$ for every $v\in V$.
In the same vein, $L-L'$ is a list assignment of $G$ such that $(L-L')(v)=L(v)\smallsetminus L'(v)$ for every $v\in V$. 
We define \emph{$\omega$ narrowed down to $L'$} as $\omega|_{L'}:V\to \mathbb{N}$ where  $\omega|_{L'}(u)= \omega(u)$ if $\omega(u)\in L'(u)$ and $\omega|_{L'}(u)= 0$ otherwise.

Suppose we are choosing a colouring $\omega$ uniformly at random from the set of all partial $L$-coloruings of $G$. We shall show that with positive probability it complies with~(\ref{MainGoal}), thus completing the proof.
For this aim we shall in particular make use of the following straightforward, yet very useful fact, 
which corresponds to the key Bernshteyn's observation utilized while adapting Molloy's proof.
\begin{eqnarray}
&&For ~any ~fixed ~L'\subseteq L ~and ~a~partial ~L'-colouring ~\omega' ~of ~G,  \nonumber\\
&& the ~partial ~colouring ~\omega|_{L-L'} ~of ~G  ~conditioned ~on ~the ~event ~\{\omega|_{L'}=\omega'\} \label{KeyObservation}\\
&&is ~uniformly ~distributed ~over ~all ~partial ~(L-L')^{\omega'}-colourings ~of ~G.\nonumber
\end{eqnarray}

For every $u\in V$ we define the following bad event:
$$A_u:~ \omega(u)=0~~ {\rm and} ~~\left(|L^\omega(u)|< (1-\varepsilon)n^\varepsilon ~~{\rm or}  ~~|N_{L^\omega}^c(u)| > \varepsilon^2n^\varepsilon ~{\rm for ~some} ~c\in L^\omega(u)\right).$$ 
Let $\mathcal{A}$ denote the set of all such events. Note that if none of these holds, then~(\ref{MainGoal}) is fulfilled, as $\varepsilon^2< (1-\varepsilon)/(2e)$ for $0<\varepsilon<1/3$. 
For every $u\in V$ define $\Gamma(A_u)=\{A_v:v\in V, {\rm dist}(v,u)\leq 3\}$, where ${\rm dist}(v,u)$ denotes the distance of $v$ and $u$ in $G$, 
hence $|\Gamma(A_u)|\leq \Delta^3$. Due to the Local Lemma, in order to complete the proof it is therefore sufficient to prove that for each $u\in V$ and any set of events $\mathcal{B}\subseteq \mathcal{A}\smallsetminus (\Gamma(A_u)\cup \{A_u\})$,
\begin{equation}\label{ConditionalProbabilityAu1}
\mathbf{Pr}\left(A_u~|~\bigcap_{A\in\mathcal{B}}\overline{A}\right) < \Delta^{-4}.
\end{equation}
Fix any $\mathcal{B}$ as above.
Note that each $A_v$ is pinned down by the colours of vertices at distance at most $2$ from $v$, hence $\mathcal{B}$ is fully determined by the colours of vertices outside $N[u]$.
It is thus sufficient to prove that the probability of $A_u$ is upper-bounded by $\Delta^{-4}$ under the condition that $\omega$ assigns any fixed set of values outside $N[u]$ -- for technical reasons we shall assume that even slightly more is prearranged. Firstly, we may assume that $\omega(u)=0$, as otherwise there is nothing to prove. 
Secondly, we shall assume some partial knowledge concerning the neighbours of $u$, namely regarding all colours outside $L(u)$ available for these.
Let $L'\subseteq L$ be the $L$-subassignment such that 
$$L'(v)=\left\{\begin{array}{lcr}
L(v) & if & v\in V\smallsetminus N(u)\\
L(v)\smallsetminus L(u) &if & v\in N(u)
\end{array}
\right..$$ 
Let $\omega'$ be any fixed partial $L'$-colouring of $G$ such that $\omega'(u)=0$. We shall prove that 
\begin{equation}\label{ConditionalProbabilityAu2}
\mathbf{Pr}\left(A_u~|~\omega|_{L'}=\omega'\right) < \Delta^{-4},
\end{equation}
which implies~(\ref{ConditionalProbabilityAu1}), thus finishing the proof.
Let  $v_1,v_2,\ldots,v_m$ be all pairwise distinct vertices in $V$ with $L_i:=(L-L')^{\omega'}(v_i)\neq \emptyset$ for every $i$.
Then $v_i\in N(u)$, and hence, $L_i\subseteq L(u)$ (by definition of $L'$) for every $i\in [m]$.
Note that as $L$ is triangle-free, there is thus no edge $v_iv_j$ with $L_i\cap L_j\neq \emptyset$ in $G$, and therefore any choice of colours in $L_i$, $i\in [m]$ is conflict-free (i.e. cannot result in the same colour associated to adjacent vertices in $G$).
Thus assuming $\omega|_{L'}=\omega'$, any vector of values in 
${^0L_1}\times{^0L_2}\times\ldots\times{^0L_m}$ can be assigned to $(v_1,v_2,\ldots,v_m)$ by 
$\omega|_{L-L'}=\omega-\omega'$ (which due to definition of $L'$ must assign $0$ to all remaining vertices), while by~(\ref{KeyObservation}), the probability of each of these to occur is uniformly distributed.  
The resulting model is thus equivalent to the new Coupon Collector setting discussed in Section~\ref{MolloysCouponCollectorsProblem}. 
Note that $|L(u)|= n = \left\lfloor(1+4\varepsilon)\Delta/\log\Delta\right\rfloor \geq (1+3\varepsilon)\Delta/\log\Delta$ (for $\Delta$ large enough), while $m\leq\Delta$, hence
\begin{eqnarray}
(1-\varepsilon)n\log n &\geq& (1-\varepsilon) (1+3\varepsilon)\frac{\Delta}{\log\Delta} (\log \Delta-\log\log\Delta)\nonumber\\
&\geq& (1-\varepsilon) (1+2\varepsilon)\Delta \geq \Delta \geq m.
\end{eqnarray}
By Lemma~\ref{CouponCollectorsWithBlanksLemma} we thus obtain that:
\begin{equation}\label{ConditionalProbabilityAu3}
\mathbf{Pr}\left(|L^\omega(u)|< (1-\varepsilon)n^\varepsilon~|~\omega|_{L'}=\omega'\right) < e^{-\frac{\varepsilon^2n^\varepsilon}{2}}
\end{equation}
and for every $c\in L(u)$:
\begin{equation}\label{ConditionalProbabilityAu4}
\mathbf{Pr}\left(c\in L^\omega(u) \wedge |N_{L^\omega}^c(u)| > \varepsilon^2n^\varepsilon~|~\omega|_{L'}=\omega'\right) < e^{-\frac{\varepsilon^2n^\varepsilon}{6}}.
\end{equation}
As $|L(u)| = n \leq \Delta$, inequality (\ref{ConditionalProbabilityAu2}) is implied by~(\ref{ConditionalProbabilityAu3}) and~(\ref{ConditionalProbabilityAu4}) for $\Delta$ large enough. 
\qed
\end{pf}

\section{Independent sets and DP-colourings}

The same result as discussed in the previous paragraph holds also within the more general environment of DP-colourings. 
In fact the presentation of the proof itself, based on Bernshteyn's idea of exploiting randomly chosen independent sets~\cite{Bernshteyn}, seems even more transparent in such a setting. 
Our extension however felt more natural within the classical list colouring environment, which urged 
us towards exposing it in details for such a setting.
Below we nevertheless argue how it can be further stretched out towards correspondence colourings. 
We shall need a few additional notations to that end. 
We first discuss how these naturally arise with reference to list colourings. 

Given a graph $G=(V,E)$ with a list assignment $L$ we define an auxiliary $L$-conflict graph $H$. 
For every $v\in V$ we denote by 
$$L_v:=\{v_c:c\in L(v)\}$$ 
the so-called \emph{$v$-set}, which includes $|L(v)|$ copies of $v$ -- each indexed by a colour from the list associated to $v$ (where $L_v=\emptyset$ if $L(v)=\emptyset$).
The vertex set of $H$ is formed of all $v$-sets, i.e. $V(H):=\bigcup_{v\in V}L_v$.
Two vertices $u_c$, $v_{c'}$ are in turn joined by an edge in $H$ either if $u=v$ or if $uv\in E$ and $c=c'$. We also denote by $H^*$ the graph obtained from $H$ by removing every edge $u_cv_{c'}$ of the first type, i.e. with $u=v$. 
Note that $L$ is a triangle-free list assignment if and only if $H^*$ has no triangles. Moreover, 
$\omega:V\to \mathbb{N}$ is a partial $L$-colouring of $G$ if and only if 
$$I_\omega:=\{v_{\omega(v)}:v\in V\smallsetminus O_\omega\}$$ 
is an independent set in $H$ (which is the more an independent set in $H^*$). 
For any independent set $\tilde{I}$ in $H$ and a subgraph $\tilde{H}$ of $H$, set:
$$\tilde{H}^{\tilde{I}}:= \tilde{H}-N_H[\tilde{I}].$$ 

Given $L'\subseteq L$ consider an $L'$-conflict graph $H'$. Then $H'$ is obviously an induced subgraph of $H$.
Let $H-H':=H-V(H')$.
Suppose $I$ is an independent set in $H$ chosen randomly and uniformly from the family of all independent sets in $H$.
We define the independent subset $I|_{H'}:=I\cap V(H')$, which is just $I$ narrowed down to $H'$.
The following straightforward observation naturally extends Bernshteyn's key observation~($\#$) from~\cite{BernshteynAsymptoticDP}:
\begin{eqnarray}
&&For ~any ~fixed ~induced ~subgraph ~H' ~of ~H ~and ~an ~independent ~set ~I' ~in ~H',   \nonumber\\
&& the ~independent ~set ~I|_{H-H'}  ~conditioned ~on ~the ~event ~\{I|_{H'}=I'\}\label{KeyObservationIndependent}\\
&&is ~uniformly ~distributed ~over ~all ~independent ~sets~in ~(H-H')^{I'}.\nonumber
\end{eqnarray}
This is actually just a variation of observation~(\ref{KeyObservation}), expressed in a new language and notation, where $I,I',I|_{H'},I|_{H-H'}$ are one-to-one twins of $\omega,\omega',\omega|_{L'},\omega|_{L-L'}$, resp., and independent sets in $(H-H')^{I'}$ correspond unambiguously to partial $(L-L')^{\omega'}$-colourings (here a vertex $v$ of $G$ coloured $0$ relates to absence of vertices from the $v$-set in the corresponding independent set). 
The proof of Theorem~\ref{Triangl-Free-Epsilon-Theorem} could thus be directly rephrased in this language. 
This setting is also useful, though not indispensable, to clearly state how our result generalizes to DP-colourings.

Note first that for any edge $uv$ in $G$ the edges between $L_u$ and $L_v$ in $H$ form a matching. This guarantees that the same colour cannot be used by adjacent vertices in $G$. Within correspondence colouring, which also refers to a given list assignment $L$, we (sometimes) admit the same colours on adjacent vertices, as a given colour $c$ from the list of a vertex $u$, if associated to $u$, excludes in this setting at most one colour, say $c'$, from the list of admissible colours for its neighbour $v$, but not necessarily colour $c$, i.e. we do not have to have $c=c'$.
(Then $c'\in L(v)$ chosen as colour for $v$ symmetrically forbids $c\in L(u)$ to be associated to $u$.) 
This translates to one slight, yet frequently consequential, change in the structure of a conflict graph $H$. Namely, for any  $uv\in E(G)$ the edges between $L_u$ and $L_v$ in $H$ must still form a matching, but an edge of such a matching does not have to join vertices indexed by the same colour.
Such $H$ (along with $L$) shall be called a cover, cf.~\cite{Bernshteyn}. 
Formally, a \emph{$k$-fold cover} of $G=(V,E)$ is a pair $\mathcal{H}=(L,H)$ where $L$ is a list $k$-assignment and $H$ is a graph with $V(H):=\bigcup_{v\in V}L_v$ such that for any $u,v\in V$ the graph $H[L_v]$ is complete while $E_H[L_u,L_v]$ forms a (possibly empty) matching if $uv\in E$ or an empty set otherwise. (The graph $H^*$ is defined as above.)
An \emph{$\mathcal{H}$-colouring} of $G$ is in turn simply defined as an independent set of size $|V|$ in $H$. 
The least $k$ such that $G$ is $\mathcal{H}$-colourable for every $k$-fold cover $\mathcal{H}$ of $G$ is called the DP-chromatic number of $G$ and denoted by $\chi_{DP}(G)$.
By definition it is at least as large as the choosability of a graph. It is however known that in general it can be far greater, as exemplify e.g. the balanced complete  bipartite graphs, for which we have  $\chi_l(K_{d,d})=(1+o(1))\log_2d$~\cite{ErdosRubinTaylor}, and $\chi_{DP}(K_{d,d})=(1/2+o(1))d/\log d$~\cite{BernshteynAsymptoticDP}.

We observe that the near-optimal upper bound on $\chi_{DP}$ for triangle-free graphs
given by Bernshteyn~\cite{Bernshteyn} (which was an improvement over his previous result from~\cite{BernshteynAsymptoticDP}) extends to the setting where no triangle is induced by the lists, while graphs themselves are allowed to contain (arbitrarily many) triangles, that is when no triangle appears in $H^*$ -- we say that the cover $\mathcal{H}$ is \emph{triangle-free} in such a case.
\begin{theorem}\label{Triangl-Free-DP-Theorem}
For any $0<\varepsilon<1/3$ there exists $\Delta_\varepsilon$ such that every graph $G$ 
with maximum degree $\Delta\geq \Delta_\varepsilon$ is $\mathcal{H}$-colourable for any triangle-free  $n$-fold cover $\mathcal{H}$ of $G$ with 
$n\geq\left\lfloor(1+4\varepsilon)\Delta/\log\Delta\right\rfloor$.
\end{theorem}
The proof of this fact is almost the same as the one of Theorem~\ref{Triangl-Free-Epsilon-Theorem}.
There are only slight adaptations to be made one needs to realize.
These are mainly three things: the new wrapping-up lemma (a correspondent of Theorem~\ref{Lists8TimesAvailiabilityObs}), the way $A_v$ translates to the new setting and an adjustment in the application of the Coupon Collector's Problem.

A more general variant of Theorem~\ref{Lists8TimesAvailiabilityObs}, suitable for our purposes takes the following form (we quote its version with a constant $1/2$, though for our goals it would also be sufficient to use its weaker version with constant $1/(2e)$, which has equally simple Local Lemma based proof as Theorem~\ref{Lists8TimesAvailiabilityObs}).
\begin{theorem}[Haxell~\cite{Haxell-ListColouring}]\label{Lists2TimesAvailiabilityLemma-Indep}
For any graph  $F$ and its vertex partition $V(F)=V_1\cup V_2\cup\ldots\cup V_t$ such that $\Delta(F)\leq l/2$ for some positive $l$ and  $|V_i|\geq l$ for $i\in[t]$, there is an independent set $I$ in $G$ with $|I\cap V_i|=1$ for $i\in [t]$.
\end{theorem}
For every $u\in V$, the event $A_u$ takes the following form in turn:
$$A_u: ~I\cap L_u=\emptyset~~ {\rm and} ~~(~|L^I_u|< (1-\varepsilon)n^\varepsilon ~~{\rm or}  
~~|N_{(H^I)^*}(u_c)| 
> \varepsilon^2n^\varepsilon ~{\rm for ~some} ~u_c\in L^I_u ~)$$ 
where 
$$L^I_u:=L_u\smallsetminus N_H[I]$$
(while as above, $H^I=H- N_H[I]$ and $(H^I)^*$ equals $H^I$ with edges inside all $L_v$'s removed).
If none of the events $A_u$ holds, the assumptions of Theorem~\ref{Lists2TimesAvailiabilityLemma-Indep} are fulfilled (for small enough $\varepsilon$) by $F=(H^I)^*$, 
and hence a desired complement of $I$ can be found.
Due to the Local Lemma it thus remains to prove that the probability of $A_u$ is small enough. 
This follows by almost the same argument as before -- it is sufficient to utilize again the fact that $H^*$ is triangle-free to translate our setting to Molloy's Coupon Collector's Problem, which requires a small additional adjustment of decks' definition. 
The vertices $v_1,v_2,\ldots,v_m\in N(u)$ and $H'$ (corresponding to $L'$) are defined analogously as previously, 
i.e. $H'=H-\bigcup_{c\in L(u)} N_{H^*}(u_c)$.
 We however  apply Lemma~\ref{CouponCollectorsWithBlanksLemma} to the following decks:
 $$L_i=\{c\in L(u):N_H(u_c)\cap V((H-H')^{I'})\cap L_{v_i}\neq \emptyset\},$$ 
 $i\in[m]$, where a given vector of choices 
 $x=(x_1,x_2,\ldots,x_m)\in {^0L_1}\times{^0L_2}\times\ldots\times{^0L_m}$ defines the independent set $I_x:=\bigcup_{i\in [m]} N_H(u_{x_i})\cap L_{v_i}$, corresponding to $I|_{H-H'}$.
The rest is essentially (or even virtually) the same as in the proof of Theorem~\ref{Triangl-Free-Epsilon-Theorem}.
We omit details here, as these introduce nothing new to the argument. 

The adjustments above might be easier to follow after reading the next section, where we adapt the new notation to prove a strengthening of the second result of~\cite{BernshteynAsymptoticDP} directly 
in the correspondence colouring setting.

\section{$K_r$-free lists with $r\geq 4$}

Generalizing the notation above, given a graph $F$, we say that a cover $\mathcal{H}=(L,H)$ of a graph $G$ is \emph{$F$-free} if the corresponding $H^*$ is $F$-free. 
Following Molloy~\cite{Molloy-small_clique_number}, in order to analyze $K_r$-free covers,
we shall utilize the following observation due to Shearer~\cite{Shearer}, where  ${\rm Ind}(F)$ denotes the family of independent sets in $F$, while ${\rm ind}(F):=|{\rm Ind}(F)|$ -- the number of these (including the empty set).

\begin{lemma}[Shearer~\cite{Shearer}]\label{ShearerLemma}
For any $r\geq 2$, if $F$ is a $K_r$-free graph, then $2^{|V(F)|}\geq {\rm ind}(F)\geq 2^{|V(F)|^\frac{1}{r-1}-1}$.
\end{lemma}
The upper bound in Lemma~\ref{ShearerLemma} refers to the number of all subsets of $V(F)$, while the lower bound has a few line straightforward inductive (with respect to $r$) proof, divided into two cases depending on whether  $F$ has a vertex of degree at least $|V(F)|^\frac{r-2}{r-1}$ or not, see~\cite{Molloy-small_clique_number,Shearer} for details.
We shall use it to show that most independent sets in a $K_r$-free graph $F$ must have relatively large cardinalities, strengthening slightly the corresponding lemmas in~\cite{Molloy-small_clique_number,BernshteynAsymptoticDP} for graphs with large enough ${\rm ind}(F)$.

\begin{lemma}\label{LemmaNoIndSets}
If $F$ is a $K_r$-free graph, then at least $(1-{\rm ind}(F)^{-\frac{1}{r^2}}){\rm ind}(F)$ independent sets in $F$ have size at least $\frac{1+\frac{1}{r}}{r-1}\log_2{\rm ind}(F)/\log_2\log_2{\rm ind}(F)$, provided that ${\rm ind}(F)$ is large enough.
\end{lemma} 

\begin{pf} Set $n:=|V(F)|$, $i:={\rm ind}(F)$ and $t:=\lfloor\frac{1+\frac{1}{r}}{r-1}\log_2i/\log_2\log_2i\rfloor$. We shall show something more, namely that the number of all subsets of $V(F)$ of cardinality at most $t$ equals no more than 
$i/i^{1/r^2}$. 
As ${n \choose j}\leq (\frac{en}{j})^j$ and since by Lemma~\ref{ShearerLemma}, $(1+\log_2i)^{r-1}\geq n \geq \log_2i\geq 2t$, for $i$ large enough,
\begin{eqnarray}
\sum_{j=0}^t{n \choose j} &\leq& (t+1)\left(\frac{en}{t}\right)^t\leq (t+1)\left(\frac{e}{t}\right)^t \left(1+\log_2i\right)^{(r-1)t} 
\leq \left(\frac{3}{t}\right)^t \left(\log_2i\right)^{(r-1)t} \nonumber\\
&\leq& \left(\frac{3(r-1)\log_2\log_2i}{\log_2i}\right)^t\left(\log_2i\right)^{(r-1)t} 
\leq\left(\frac{1}{\left(\log_2i\right)^\frac{r-1}{r}}\right)^t \left(\log_2i\right)^{(r-1)t} \nonumber\\
&=& \left(\log_2i\right)^{(1-\frac{1}{r})(r-1)t}
\leq  \left(\log_2i\right)^{(1-\frac{1}{r})(1+\frac{1}{r})\frac{\log_2i}{\log_2\log_2i}} = i^{1-\frac{1}{r^2}}, \nonumber
\end{eqnarray}
as claimed.
\qed
\end{pf}

We now extend the result of Bernshteyn~\cite{BernshteynAsymptoticDP}  towards $K_r$-free covers (rather than just $K_r$-free graphs), improving at the same time the multiplicative constant in the corresponding upper bound. 
One may note it could still be slightly improved with this technique,
 but we do not strive to do this,
 as its order might still not be optimal, cf.~\cite{AlonKrivelevichSudakov}.
\begin{theorem}\label{Kr-Free-DP-Theorem}
For every $r\geq 4$ there exists a constant $\Delta_r$ such that each graph $G$ 
with maximum degree $\Delta\geq \Delta_r$ is $\mathcal{H}$-colourable for any
$K_r$-free  $k$-fold cover  $\mathcal{H}$ of $G$ with 
$k\geq 2(r-2)\Delta\frac{\log_2\log_2\Delta}{\log_2\Delta}$.
\end{theorem}

\begin{pf}
Let $G=(V,E)$ be a graph of sufficiently large maximum degree $\Delta$ with respect to a fixed $r\geq 4$, i.e. large enough so that all inequalities below hold.
Let $\mathcal{H}=(L,H)$ be a $K_r$-free $k$-fold cover of $G$ with $$k\geq 2(r-2)\Delta\frac{\log_2\log_2\Delta}{\log_2\Delta}.$$ 
Suppose we are choosing $I$ independently at random from the family of all independent sets in $H$; denote by $O_I$ the set of all vertices $v\in V$ with $L_v\cap I=\emptyset$.
Set 
$G_I:=G[O_I]$ and
$$l:=\Delta^{\frac{1}{2}+\frac{1}{8r}}.$$ 
We shall show that with positive probability, for every $u\in O_I$: 
\begin{equation}\label{MainGoodEvent}
|L^I_u|\geq l~~~~ {\rm and}~~~~ d_{G_I}(u)<l. 
\end{equation}
This is sufficient, as one may then finish up the construction of an $\mathcal{H}$-colouring, i.e. $|V|$-element independent set in $H$, greedily adding one by one available elements from all remaining non-empty $v$-sets.
We denote the following bad events for every $u\in V$:
\begin{eqnarray}
B^1_u: && u\in O_I~~ {\rm and} ~~|L^I_u|< l ; \nonumber\\ 
B^2_u: && d_{G_I}(u)\geq l~~ {\rm and} ~~|L^I_v|\geq l
~{\rm for } ~v\in N_{G_I}(u). \nonumber
\end{eqnarray}
Denote by $\mathcal{B}$ the set of all such events.  
Note that if none of these holds, then~(\ref{MainGoodEvent}) is fulfilled for each  $u\in O_I$.
Let for $\alpha=1,2$, $\Gamma(B^{\alpha}_u):=\{B^1_v,B^2_v:v\in V, {\rm dist}(v,u)\leq 3\}$. 
Due to the Local Lemma, in order to prove the theorem it thus suffices 
to show that for every $u\in V$, $\alpha\in\{1,2\}$ and any set of events $\mathcal{A}\subseteq \mathcal{B}\smallsetminus (\Gamma(B^{\alpha}_u)\cup \{B^{\alpha}_u\})$:
\begin{equation}\label{ProbBuGoal}
\mathbf{Pr}\left(B^{\alpha}_u~|~\bigcap_{B\in\mathcal{A}}\overline{B}\right) < \Delta^{-4}.
\end{equation}

We consider $\alpha=1$ first.
Let us fix any $u\in V$ and let analogously as before $L'\subseteq L$ be defined as follows: 
$$L'(v)=\left\{\begin{array}{lcr}
L(v) & if & v\in V\smallsetminus N(u)\\
\{c\in L(v):N_H(v_c)\cap L_u=\emptyset\} 
&if & v\in N(u)
\end{array}
\right..$$
Let $\mathcal{H}'=(L',H')$ be a subcover of $\mathcal{H}$ induced by $L'$, i.e. $H'=H[\bigcup_{v\in V}L'_v]$.  
Since events in any $\mathcal{A}\subseteq \mathcal{B}\smallsetminus (\Gamma(B^1_u)\cup \{B^1_u\})$ are determined by $I\cap \bigcup_{v\in V\smallsetminus N[u]}L_v$ while the probability of $B^1_u$ is zero if $L_u\cap I\neq \emptyset$, in order to show~(\ref{ProbBuGoal}) for $\alpha=1$, it is sufficient to prove that
\begin{equation}\label{ProbBuGoal2}
\mathbf{Pr}\left(B^1_u~|~I|_{H'}=I'\right) < \Delta^{-4}
\end{equation}
for any fixed independent set $I'$ in $H'$ (thus also independent in $H$, as $H'$ is an induced subgraph of $H$) such that $L_u\cap I'=\emptyset$. 

Set $H_0:=H-H'$. 
By~(\ref{KeyObservationIndependent}), $I|_{H_0}$ is uniformly distributed over all independent sets in $H_0^{I'}$ (conditioned on the event: $\{I|_{H'}=I'\}$). 
Note that (given a fixed $I'$) $I|_{H_0}$ may be chosen via the following random process, 
where we shall first choose it uniformly at random, and then subsequently resample our choice in the neighbourhoods of the consecutive elements in $L_u$.
Let $L(u)=\{c_1,c_2,\ldots,c_k\}$ and set $N_i:=N_{H^*}(u_{c_i})\cap V(H_0^{I'})$ for $i=1,\ldots,k$. Let $I_0$ be chosen uniformly at random from the family of all independent sets in $H_0^{I'}$. For $i=1,\ldots,k$ we then proceed as follows.
\begin{itemize}
\item[(a)] Let $H_i:=(H_0[N_i])^{I_{i-1}\smallsetminus N_i}$; note $H_i$ is in particular a $K_{r-1}$-free  induced subgraph of $H_0^{I'}$ (if there was a clique of order $r-1$ in $H_i$, then its vertex set along  
with $u_{c_i}$ would induce $K_r$ in $H^*$); 
\item[(b)] Choose uniformly at random an independent set $I'_i$ in $H_i$;
\item[(c)] Set $I_i:=(I_{i-1}\smallsetminus N_i)\cup I'_i$.
\end{itemize}
By~(\ref{KeyObservationIndependent}), it is clear that the resulting $I_k$ is uniformly distributed over the independent sets in $H_0^{I'}$.
Note that the choices performed in (b) within the process above can be realized via the following randomized procedure.
Prior to launching the whole resampling algorithm we evaluate $2k$ independent random variables with uniform distribution over $[0,1)$ interval: $X_1,\ldots,X_k,Y_1,\ldots,Y_k\sim U[0,1)$. We may 
assume that there is a fixed ordering $\preceq$ in any family of independent sets in $H$, arranging these in size-nondecreasing manner. Let $k_x$ be an auxiliary variable counting how many times we shall be using variables $X_j$ (this shall refer to the number of times when ${\rm ind}(H_i)$ is relatively small), and we initially set $k_x:=0$. Suppose that in a given consecutive step $i\in[k]$, ${\rm Ind}(H_i)=\{\tilde{I}_p\}_{1\leq p\leq l_i}$ where $\emptyset=\tilde{I}_1\preceq \tilde{I}_2\preceq\ldots\preceq \tilde{I}_{l_i}$. We then proceed as follows: 
\begin{itemize}
\item[(b1)] if $l_i< \Delta^{\frac{1}{2}-\frac{1}{4r}}$, then we increase $k_x$ by $1$ and 
exploit $X_{k_x}$ to establish $I\cap N_i$, setting $I'_i=\tilde{I}_q$ iff $X_{k_x}\in [\frac{q-1}{l_i},\frac{q}{l_i})$; 
note that if $X_{k_x}<\Delta^{-\frac{1}{2}+\frac{1}{4r}}$, then $I'_i=\emptyset$;
\item[(b2)] if $l_i\geq \Delta^{\frac{1}{2}-\frac{1}{4r}}$, 
we exploit $Y_{i-k_x}$ in turn, and set 
$I'_i=\tilde{I}_q$ iff $Y_{i-k_x}\in [\frac{q-1}{l_i},\frac{q}{l_i})$; 
note that if $Y_{i-k_x}\geq (\Delta^{\frac{1}{2}-\frac{1}{4r}})^{-\frac{1}{(r-1)^2}}$, then by Lemma~\ref{LemmaNoIndSets}, due to the fact that $H_i$ is $K_{r-1}$-free,
$|I'_i|\geq \frac{1+\frac{1}{r-1}}{r-2}\log_2l_i/\log_2\log_2l_i 
\geq \frac{(1+\frac{1}{r-1})(\frac{1}{2}-\frac{1}{4r})}{r-2}\log_2\Delta/\log_2\log_2\Delta 
= \frac{1+\frac{1}{2(r-1)}}{2(r-2)}\log_2\Delta/\log_2\log_2\Delta$.
\end{itemize}
Let from now on $k_x$ refer to the value of this parameter at the end of the process. Set 
\begin{eqnarray}
p_1&:=& \Delta^{-\frac{1}{2}+\frac{1}{4r}}, \nonumber\\ 
p_2&:=&1-\Delta^{-r^{-3}}<1-(\Delta^{\frac{1}{2}-\frac{1}{4r}})^{-\frac{1}{(r-1)^2}}, \nonumber\\ 
k_1&:=& \lceil 2\Delta^{1-\frac{1}{8r}} \rceil, \nonumber\\  
k_2&:=& k-k_1\geq k(1-\Delta^{-r^{-3}})=kp_2. \nonumber 
\end{eqnarray}
Let $X$ be the number of $X_j$'s with $1\leq j\leq k_1$ such that 
$X_j<p_1$ and let $Y$ be the number of $Y_j$'s with $1\leq j\leq k_2$ such that 
$Y_j\geq 1-p_2$. 
Note that by the Chernoff Bound:
\begin{eqnarray}
\mathbf{Pr}(X<l \vee Y<p_2^3k) &\leq& \mathbf{Pr}(X<0.5k_1p_1) + \mathbf{Pr}(Y<p_2^2k_2) \nonumber\\
&\leq& e^{-\frac{k_1p_1}{8}} + e^{-\frac{k_2p_2}{2\Delta^{2r^{-3}}}} ~<~ \Delta^{-4}. \label{XYprobability}
\end{eqnarray}
Thus with probability at least $1-\Delta^{-4}$, $X\geq l$ and $Y\geq p_2^3k$.
Suppose both of these inequalities hold.
If $k_x<k_1$, then the values of $Y_1,\ldots,Y_{k_2}$ were utilized within (b2), and hence, 
by the definition of $Y$ and 
the fact that $1-p_2 > (\Delta^{\frac{1}{2}-\frac{1}{4r}})^{-\frac{1}{(r-1)^2}}$, we have $|I'_i|\geq  \frac{1+\frac{1}{2(r-1)}}{2(r-2)}\log_2\Delta/\log_2\log_2\Delta$ for at least $p_2^3k$ values of $i$. Thus:
\begin{eqnarray}
\Delta &=& |N(u)|\geq |I\cap(\bigcup_{v\in N(u)}L_v)| 
\geq \sum_{i=1}^k|I'_i| \geq p_2^3k\cdot \frac{1+\frac{1}{2(r-1)}}{2(r-2)}\frac{\log_2\Delta}{\log_2\log_2\Delta} \nonumber\\
&\geq& p_2^32(r-2)\Delta\frac{\log_2\log_2\Delta}{\log_2\Delta}\frac{1+\frac{1}{2(r-1)}}{2(r-2)}\frac{\log_2\Delta}{\log_2\log_2\Delta}
= p_2^3 \left(1+\frac{1}{2(r-1)}\right) \Delta >\Delta \nonumber
\end{eqnarray}
for $\Delta$ large enough, which is a contradiction.
 It follows that $k_x\geq k_1$. Then by (b1) and the definition of $X$, $I'_i=\emptyset$ for at least $l$ values of $i$, and hence $|L^I_u|\geq l$. 
Thus~(\ref{XYprobability}) implies~(\ref{ProbBuGoal2}). 

It is left to discuss $\alpha=2$. Consider $u\in V$ and any set $\{v_1,\ldots,v_{\lceil l \rceil}\}$ of $\lceil l\rceil$ neighbours of $u$ in $G$. 
In order to prove~(\ref{ProbBuGoal}) it is sufficient to show that for any fixed independent set $I''$ in 
$H$ which is disjoint with $L_{v_j}$ for $j=1,\ldots \lceil l \rceil$:
 \begin{equation}\label{ProbBuPart2Spec}
\mathbf{Pr}\left(|L^I_{v_j}
|\geq l~{\rm for } ~j=1,\ldots,\lceil l\rceil~|~I\smallsetminus\bigcup_{j=1}^{\lceil l\rceil}L_{v_j}=I''\right) < {\Delta\choose {\lceil l\rceil}}^{-1}\Delta^{-4}.
\end{equation}
In order to prove this in turn, note that if  $|L^{I''}_{v_j}|<l$ for some $j\in\{1,\ldots,\lceil l\rceil\}$, then the probability above is $0$. Otherwise, as the event 
$\{|L^I_{v_j}|\geq l~{\rm for } ~j=1,\ldots,\lceil l\rceil\}$ 
implies that $I\cap\bigcup_{j=1}^{\lceil l\rceil} L^{I''}_{v_j} 
=\emptyset$, while there are at least $\lceil l\rceil !$ different independent sets in $H[\bigcup_{j=1}^{\lceil l\rceil} L^{I''}_{v_j} 
]$,
then Stirling's formula combined with inequality ${n \choose m} \leq (\frac{en}{m})^m$ imply that:
\begin{eqnarray}
&&\mathbf{Pr}\left(|L^I_{v_j}| 
\geq l~{\rm for } ~j=1,\ldots,\lceil l\rceil~|~I\smallsetminus\bigcup_{j=1}^{\lceil l\rceil}L_{v_j}=I''\right) \nonumber\\
&\leq& \frac{1}{\lceil l\rceil !}  
<\frac{1}{\left(\frac{\lceil l\rceil}{e}\right)^{\lceil l\rceil}}
<\frac{1}{\left(\frac{e\Delta}{\lceil l\rceil}\right)^{\lceil l\rceil}\Delta^4}  
\leq \frac{1}{{\Delta\choose {\lceil l\rceil}}\Delta^{4}}
\end{eqnarray}
for large enough $\Delta$, as desired.
\qed
\end{pf}

\end{document}